\documentstyle[12pt]{article}
\title{\bf On Integrable Solutions of
Impulsive Delay Differential Equations}

\author{
L. Berezansky $^{1}$
\\Ben-Gurion University of the Negev, \\
Department of Mathematics and Computer Science, \\
Beer-Sheva 84105, Israel, \\
E. Braverman $^{2}$ \\
Technion - Israel Institute of Technology,\\
Department of Mathematics, 32000, Haifa, Israel }

\begin{document}
\maketitle

\footnotetext[2]{Supported by : The Centre for
 Absorption in Science, Ministry of Immigrant Absorption State of Israel }
\footnotetext[1]{Supported by Israel Ministry of Schience and
Technology}

\section{Introduction}

{}~~~~~The belonging of solutions to a certain function space
is a characteristic property for studying the asymptotic
behavior of solutions of differential equations.
Many works are concerned with the connection between the properties
of solutions and stability.
We name here the monographs [1-3] on ordinary differential equations
and the works [4-9] on functional differential equations.
For differential equations with impulses this problem was investigated
in [10-12] for ordinary differential equations and in [13] for equations
with delay.

The present paper deals with the following problems:

admissibility of a pair of spaces for a differential operator,
i.e. action conditions for this operator in corresponding
function spaces;

admissibility of a pair of spaces for a differential equation,
i.e. the conditions of belonging of all solutions to a certain space
if provided that the right hand side belongs to the other space;

connection between admissibility and exponential stability
for impulsive differential equations.

All function spaces considered are the space of locally integrable functions
and its subspaces.
Explicit conditions for existence of integrable solutions and
for exponential stability are obtained as corollaries of these results.

The present paper is organized as follows.
In section 2 the equation studied is described and the hypotheses
are introduced.
Section 3 deals with auxiliary results.
In particular the solution representation formula is given and
the properties of certain spaces of differentiable on the half-line
functions are described.
The proofs of these results are presented in the last section 7.
In section 4 admissibility of a pair of spaces is considered.
In section 5 stability problems are investigated.
Finally, section 6 gives explicit stability results.

In conclusion we outline that the present work can be
treated as  [13] continued.
This paper dealt with the same problems in
the space of essentially bounded functions.

\section{Preliminaries}

{}~~~~~Let $0 = \tau_0 < \tau_1 < \dots $ be the fixed points,
$\lim _{j \rightarrow \infty} \tau _j = \infty , $
${\bf R}^n$ be the space of $n$-dimensional column vectors
$x = col (x_1, \dots ,x_n)$ with the norm
$\parallel x \parallel = \max _{1 \leq i \leq n} \mid x_i \mid $,
by the same symbol $\parallel \cdot \parallel$ we
denote the corresponding matrix norm,

$E_n$ is an $n \times n$ unit matrix,

$\chi _e : [0, \infty ) \rightarrow {\bf R}
$ is the characteristic function of the set $e :
\chi_e (t) = 1$, if $t \in e, $ and $\chi_e (t) = 0, $
otherwise.

$\bf L$ is a space of Lebesgue measurable functions
$x: [0, \infty) \rightarrow {\bf R}^n$ integrable
on any finite segment $[t, t+1]$,

${\bf L}_{\infty} \subset {\bf L}~ $ is a Banach space of essentially bounded
functions $x: [0, \infty ) \rightarrow
{\bf R}^n , \parallel \! x \! \parallel _{{\bf L}_{\infty}} =
vraisup _ {t \geq 0} \parallel x(t) \parallel ,$

\newcommand{\lp}{{\bf L}_p}
\newcommand{\depe}{{\bf D}_p}
\newcommand{\mepe}{{\bf M}_p}

$\lp  \subset {\bf L} ~(1 \leq p < \infty )$
is a Banach space of functions $x: [0, \infty) \rightarrow
{\bf R}^n$ such that $\int_0^{\infty} \parallel x(t) \parallel ^p dt
< \infty $, with a norm
$$ \parallel x \parallel _{{\bf L}_p} =
\left( \int_0^{\infty} \parallel x(t) \parallel ^p dt \right) ^{1/p}, $$

$\mepe \subset {\bf L}$ is a Banach space of functions
$x: [0, \infty) \rightarrow {\bf R}^n$ such that
$$\mu = \sup_{t>0} \left( \int_t^{t+1}
\parallel x(t) \parallel ^p dt \right) ^{1/p} < \infty, ~
1 \leq p < \infty, ~ \parallel x \parallel _{\mepe} = \mu. $$

${\bf PAC} (\tau _1, \dots , \tau_j, \dots ) $
is a linear space of
functions $x: [0, \infty) \rightarrow {\bf R}^n$
absolutely continuous on each interval $[\tau_j, \tau_{j+1} )$,
with jumps at the points $\tau_j$.
We assume that functions in {\bf PAC} are right continuous.

The same function spaces will be considered for intervals different
from $[0, \infty)$ if it does not lead to misunderstanding.

For spaces of matrix valued functions we use the same notation as for
vector valued functions.

We consider a  delay differential equation
\begin{equation}
\dot{x} (t) + \sum_{k=1}^m {A_k (t) x[h_k(t)]} =
f(t),~ t \geq 0, ~ x(t) \in {\bf R}^n , \label{e1}
\end{equation}
\begin{equation}
x(\xi ) = \varphi (\xi), \xi < 0, \label{e2}
\end{equation}
with impulsive conditions
\begin{equation}
x(\tau _j) = B_j x(\tau _j - 0) + \alpha_j ,~ j = 1,2, \dots,  \label{e3}
\end{equation}
under the
following assumptions:

(a1) $0 = \tau_0 < \tau_1 < \tau _2 < \dots $ are fixed points,
$\lim_{j \rightarrow \infty} \tau _j  = \infty $ ;

(a2) $ f \in {\bf L}, ~ A_k \in {\bf L}, ~k=1,2, \dots, m $ ;

(a3) $h_k : [0,\infty) \rightarrow {\bf R} $
are Lebesgue measurable functions, $$h_k (t) \leq t,~
k = 1, \dots , m;$$

(a4) $\varphi : (- \infty, 0) \rightarrow {\bf R}^n $
is a Borel measurable bounded function;

(a5) $ B_j \in {\bf R}^{n \times n}, ~B = \sup_j \parallel B_j \parallel
< \infty;$

(a6) $K =  \sup _{t,s >0}
\left\{ \frac{i(t,s)}{t-s}, ~i(t,s) \neq 1 \right \} < \infty. $

Here $i(t,s)$ is a number of points $\tau _j$ belonging
to the interval $(s,t).$

We denote $b = \max \{B, 1 \}, ~ I = max \{K ,1 \}. $
\vspace{5 mm}

\underline{\sl Remark.}
One can easily see that (a6) is satisfied if
$\tau_{j+1} - \tau_j \geq \rho > 0$.

\underline{{\sl Definition .}}
 A function $x \in {\bf PAC}$
is said to be {\bf a solution of the impulsive equation
(1),(2),(3)}  with the initial function $\varphi (t)$
if (1) is satisfied for almost all $t \in [0, \infty)$
and the equalities (3) hold.
\vspace{5 mm}

Below we use a linear differential operator
\begin{equation}
({\cal L} x)(t) = \dot{x} (t) + \sum_{k=1}^m {A_k (t) x[h_k(t)]} ,
{}~x(\xi) = 0, ~ \xi < 0. \label{e4}
\end{equation}

\section{Auxiliary results}

{}~~~~~In  [13] the solution representation formula for (\ref{e1})-(\ref{e3})
is presented if provided that more restrictive conditions than (a1)-(a6)
hold.
Precisely, instead of (a2) it was assumed that $f$ and $A_k$ are in
${\bf L}_{\infty}$.
However the proof of this formula preserves in the more
general case $f,A_k \in {\bf L}$.
Thus the following result is valid.

\newtheorem{uess}{Lemma}
\begin{uess}
{\em [13]} ~Suppose the hypotheses (a1)-(a6) hold.

Then there exists one and only one solution of the equation (\ref{e1})
-(\ref{e3})
satisfying $x(0)= \alpha_0$ and it can be presented as
\begin{equation}
x(t) = \int_0^t {X(t,s) f(s) ds}
- \sum_{k=1}^m {\int_0^t {X(t,s) A_k(s) \varphi [h_k(s)] ds}} +
\sum_{0 \leq \tau_j \leq t} {X(t, \tau_j) \alpha _j}, \label{e5}
\end{equation}
with $\varphi (\zeta) = 0,$ if $\zeta \geq 0. $

The matrix $X(t,s)$ in (\ref{e5}) for a fixed $s$ as a function of $t$
is a solution of the problem
$$
\dot{x} (t) + \sum_{k=1}^m {A_k (t) x[h_k(t)]} =
0,~ t \geq s, ~ x(t) \in {\bf R}^{n \times n},
$$  $$
x(\xi ) = 0,~ \xi < s,~ x(s)=E_n; ~
x(\tau _j) = B_j x(\tau _j - 0), ~ \tau_j > s    .
$$
We assume $X(t,s) = 0, ~t<s$.
\end{uess}

\underline{Definition.}
The matrix $X(t,s)$ is said to be {\bf a fundamental matrix},
$X(t,0)$ is said to be {\bf a fundamental solution}.
An operator
$$(Cf)(t) = \int_0^t X(t,s) f(s) ds $$
is said to be {\bf a Cauchy operator} of the equation (1)-(3).

For studying the equation (\ref{e1})-
(\ref{e3}) we introduce an auxiliary equation
\begin{equation}
({\cal L}_0 x)(t) \equiv
\dot{x} (t) + a x(t) = z(t), ~t \geq 0, ~x(t) \in {\bf R}^n,
\label{e6}
\end{equation}
\begin{equation}
x(\tau_j) = B_j x(\tau_j - 0). \label{e7}
\end{equation}
By
$$(C_0 z)(t) = \int_0^t X_0 (t,s) z(s) ds $$
the Cauchy operator of the equation (\ref{e6}),(\ref{e7}) is denoted.

\begin{uess}
{\em [13]} ~Suppose the hypotheses (a5) and (a6) hold
and $\nu =  a- I \ln b > 0$.

Then
$$ \parallel X_0 (t,s) \parallel \leq e^{- \nu (t-s)}. $$
\end{uess}

For each space $\lp$ we construct a subspace of ${\bf PAC}$
as follows.
Denote by $\depe$ a linear space of functions $x \in {\bf PAC}$
satisfying (\ref{e7}) and such that $x \in \lp, ~\dot{x} \in \lp $.
This space is normed, with a norm
$$ \parallel x \parallel _{\depe} = \parallel x \parallel _{\lp} +
\parallel \dot{x} \parallel _{\lp} . $$

\begin{uess}
Suppose the hypotheses (a5) and (a6) hold.

Then $\depe, ~1 \leq p < \infty,$ is a Banach space.
\end{uess}

The proof is presented in section 7.

\underline{\sl Remark.} Lemma 3 remains valid
if $\lp$ is changed by a Banach space ${\bf B} \subset {\bf L}$
if provided that the topology in ${\bf B}$ is stronger than the
topology in ${\bf L}$.
In particular ${\bf B} = {\bf L}_{\infty}$ or
${\bf B} = {\bf M}_p$ are suitable.

The following assertion supplements Lemma 3.

\begin{uess}
Suppose the hypotheses (a5) and (a6) hold and $a- I \ln b > 0$.

Then the set
$\tilde{\depe} = \{ x \in {\bf PAC} \mid ~ \dot{x} + ax \in \lp ,
{}~x(\tau_j)=B_j x(\tau_j -0) \} $
coincides with $\depe$, and the norm
\begin{equation}
\parallel x \parallel _{\tilde{\depe}} = \parallel x(0) \parallel +
\parallel \dot{x} + ax \parallel _{\lp}   \label{e8}
\end{equation}
is equivalent to the norm $\parallel \cdot \parallel _{\depe} $.
\end{uess}

The proof is also in section 7.

\section{Admissibility of pairs}

{}~~~~~\underline{\sl Definition.}
The pair ($\depe, \lp$) is said to be {\bf admissible
for a differential operator } ${\cal L}: {\bf PAC} \rightarrow
{\bf L}$ if ${\cal L} (\depe) \subset \lp$.

\underline{\sl Definition.}
Suppose the initial function $\varphi$ satisfies the hypothesis (a4)
and it is fixed.
The pair $(\lp, \depe)$ is said to be {\bf admissible} for the equation
(1)-(3) if for any $f \in \lp, ~ \alpha_j \in {\bf R}^n$
the solution is in $\depe$.

The pair $(\lp, \depe)$ is said to be
{\bf admissible on the whole} for the equation
(\ref{e1})-(\ref{e3})
if for any $f \in \lp, ~ \alpha_j \in {\bf R}^n$
and any initial function  $\varphi$
satisfying (a4) the solution is in $\depe$.

\underline{\sl Remarks.}
1. For ordinary differential equations the admissibility of the pair
$(\lp, {\bf L}_{\infty})$ is usually considered.
However this admissibility is the consequence of the admissibility of  pair
$(\lp, \depe)$.
In fact if $x \in \depe$ then for $a \in {\bf R}
 ~~~~\dot{x}~+~ax~\in~\lp$.
Let $a-I \ln~b~>~0$.
Then by Lemma 2~ $x~\in~{\bf L}_{\infty}$,
therefore the pair $(\lp,{\bf L}_{\infty})$ is admissible
for the differential equation.
Besides this, under our approach admissibility of the pair $(\lp,
\depe )$
is treated more naturally than of the pair $(\lp, {\bf L}_{\infty})$.

2. It is to be noted that the recent monograph of C.Corduneanu
[9] deals
with admissibility of pairs of spaces for integrodifferential equations
(and for general functional differential equations as well).

\vspace{5 mm}

Consider operators
\begin{equation}
(Hx)(t) = \sum_{k=1}^m A_k(t) x[h_k(t)];
{}~x(\xi)=0,~\xi <0, \label{e9}
\end{equation}
$$({\cal L} x)(t) = \dot{x}(t)+(Hx)(t) .$$

Under the hypotheses (a1)-(a3), (a5)-(a6) $H$ acts from
${\bf PAC}$ to ${\bf L}$.

\newtheorem{guess}{Theorem}
\begin{guess}
Suppose the hypotheses (a1)-(a3), (a5),(a6) hold and there exists
$\nu > 0$ such that $A_k^{\nu} \in \mepe$, where
$$A_k^{\nu} (t) = e^{\nu [t-h_k (t)]} A_k (t), ~1 \leq p< \infty
. $$

Then operators $H$ and ${\cal L}$ act from $\depe$ to $\lp$
and they are bounded.
\end{guess}

{\sl Proof.}
Let $a = \nu + I \ln b $ and $x \in \depe $ .
Then $z = \dot{x} + ax \in \lp$ and $x$ can be presented as
$$x(t)=X_0(t,0)x(0)  + \int_0^t X_0(t,s)z(s)ds . $$

In sequel $y(h(t))=0$, if $h(t)<0$,
and $a^+ = \max \{ a,0 \} $.

Thus we obtain
$$ (Hx)(t) = $$
\begin{equation}
= \sum_{k=1}^m A_k(t) X_0(h_k(t),0) x(0) +
\sum_{k=1}^m \int_0^{h_k^+ (t)} \!\!\!\! A_k(t) X_0(h_k(t),s) z(s)ds .
\label{e10}
\end{equation}

First we will obtain that a matrix valued function
$$F(t) = \sum_{k=1}^m A_k (t) X_0(h_k(t),0) $$
is in $\lp$.
To this end by Lemma 2
$$ \parallel A_k (t) X_0 (h_k(t),0) \parallel \leq
\parallel A_k (t) \parallel e^{- \nu h_k (t)} = $$
$$ = \parallel A_k (t) e^{\nu (t-h_k(t))} \parallel
e^{- \nu t} = \parallel A_k^{\nu } (t) \parallel e^{- \nu t} . $$

Therefore
$$ \int_0^{\infty} \parallel A_k^{\nu} (t) \parallel ^p e^{- \nu pt}
dt \leq \sup_{n \geq 0} \int_n^{n+1} \parallel A_k^{\nu} (t) \parallel
^p dt \sum_{n=0}^{\infty} e^{- \nu pn} \leq $$
$$ \leq \frac{\parallel A_k^{\nu} \parallel_{{\bf M}_p}^p }
{1- e^{- \nu p}} . $$

Hence $F \in \lp$.

Denote
$$ (Pz)(t) = \sum_{k=1}^m
\int_0^{h_k^+ (t)} A_k (t) X_0
(h_k(t),s) z(s) ds. $$

We will prove that $P$ acts in $\lp$ and it is bounded.
To this end
$$\parallel (Pz)(t) \parallel \leq \sum_{k=1}^m
\int_0^{h_k^+ (t)} \!\!\! \!\! \parallel A_k (t) e^{\nu [t- h_k (t)]}
\parallel
{}~ e^{- \nu (t-s)} \parallel z(s) \parallel ds = $$
$$ = \sum_{k=1}^m \int_0^{h_k^+ (t)}
\parallel A_k^{\nu} (t) \parallel ~e^{- \nu (t-s)}
\parallel z(s) \parallel ~ds. $$

Let $p=1$.
Then
$$  \parallel Pz \parallel _{{\bf L}_1}  \leq
\sum_{k=1}^m \int_0^{\infty} \int_0^t
\parallel A_k^{\nu} (t) \parallel
{}~e^{- \nu (t-s)} \parallel z(s) \parallel ~ds~dt = $$
$$ = \sum_{k=1}^m \int_0^{\infty} \left(
\int_s^{\infty} \parallel A_k^{\nu} (t) \parallel ~e^{- \nu (t-s)}
dt \right) \parallel z(s) \parallel ~ds. $$

Since
$$ \int_s^{\infty} \parallel A_k^{\nu} (t) \parallel
{}~e^{- \nu (t-s)} dt \leq
\sum_{n=[s]}^{\infty} \int_n^{n+1}
\parallel A_k^{\nu} (t) \parallel ~e^{- \nu (t-s)} dt \leq $$
$$ \leq e^{\nu s} \sum_{n=[s]}^{\infty}
e^{- \nu n} \int_n^{n+1} \parallel A_k^{\nu} (t) \parallel ~dt
\leq
e^{\nu s} \parallel A_k^{\nu} \parallel_{{\bf M}_1}
\sum_{n=[s]}^{\infty} e^{- \nu n} \leq  $$
$$ \leq \parallel A_k^{\nu} \parallel_{{\bf M}_1} \frac{e^{\nu}}
{1- e^{- \nu}}, $$
then
$$\parallel Pz \parallel_{{\bf L}_1} \leq
\frac{e^{\nu}}{1-e^{- \nu}}
\sum_{k=1}^m \parallel A_k^{\nu} \parallel_{{\bf M}_1}
\parallel z \parallel_{{\bf L}_1}. $$
Here $[s]$ is the greatest integer not exceeding $s$.

Let $1<p< \infty $.
Then similarly we obtain
$$ \parallel Pz \parallel_{\lp} \leq
 \sum_{k=1}^m  \left[ \int_0^{\infty}
\parallel A_k^{\nu } (t) \parallel^p
\left( \int_0^t e^{- \nu (t-s)} \parallel z(s) \parallel ~ds
\right)^p dt\right] ^{1/p} = $$
$$ =  \sum_{k=1}^m \left[ \int_0^{\infty} \parallel A_k^{\nu} (t)
\parallel^p \left( \int_0^t e^{- \nu (t-s)/2}
e^{- \nu (t-s)/2} \parallel z(s) \parallel ~ds \right) ^p dt
\right] ^{1/p} \leq $$
$$ \leq  \sum_{k=1}^m \left[ \int_0^{\infty}
\parallel A_k^{\nu} (t) \parallel ^p
\left( \int_0^{\infty} e^{- \nu q (t-s)/2} ds \right) ^{p/q}
\right. \times $$
$$ \times
\left. \left( \int_0^t e^{- \nu p (t-s)/2} \parallel z(s) \parallel ^p ~ds
\right) ~dt
\right] ^{1/p} \leq $$
$$ \leq \left( \frac{2}{\nu q} \right) ^{1/q}
\sum_{k=1}^m \left[ \int_0^{\infty} \int_s^{\infty}
\parallel A_k^{\nu} (t) \parallel ^p
e^{- \nu p (t-s)/2}
\parallel z(s) \parallel ^p dt~ds \right] ^{1/p}, $$
where $q = p/(p-1) $.

By repeating the previous argument we obtain
$$ \parallel Pz \parallel
_{\lp} \leq
\left( \frac{2}{\nu q} \right) ^{1/q}
\frac{e^{\nu /2}}{(1 - e^{- \nu p/2})^{1/p}}
\sum_{k=1}^m \parallel A_k^{\nu} \parallel_{\mepe}
\parallel z \parallel_{\lp} . $$
Therefore $Pz \in \lp$ and operator $P: \lp \rightarrow \lp$
is bounded.

Operator $H$ defined by (\ref{e9}) in view of (\ref{e10}) can be presented
as
$$(Hx)(t) = F(t)x(0) +(Pz)(t),~\mbox{~where~} z=\dot{x} + ax . $$
Since
$$\parallel Hx \parallel_{\lp} \leq
\parallel F \parallel_{\lp} \parallel x(0) \parallel +
\parallel P \parallel _{\lp \rightarrow \lp}
\parallel \dot{x} + ax \parallel _{\lp} \leq $$
$$ \leq \max \{
 \parallel F \parallel_{\lp},
\parallel P \parallel_{\lp \rightarrow \lp} \} \parallel x
\parallel_{\tilde{\depe} }, $$
then by Lemma 4 $H$ acts from $\depe$ to $\lp$
and it is bounded.
One can easily see that the admissibility of the pair $(\depe,
\lp)$  for the operator ${\cal L}$ is equivalent to admissibility
of this pair for $H$.
The proof of the theorem is complete.
\vspace{5 mm}

\underline{\sl Corollary.}
Suppose the hypotheses (a1)-(a3), (a5),(a6) hold,
$A_k \in {\bf M}_p, ~1 \leq p < \infty$ and there exists $\delta
>0$ such that $t-h_k(t) < \delta,~k=1, \dots, m$.

Then $H$ acts from $\depe$ to $\lp$ and it is bounded.

\vspace{5 mm}

Now we proceed to $(\lp, \depe)$ admissibility conditions
for the problem (\ref{e1}) - (\ref{e3}).
To this end consider an auxiliary equation of the type (\ref{e1}),
(\ref{e2})
$$\dot{x}(t) + \sum_{k=1}^r H_k (t) x[g_k(t)] = f(t),
{}~t \geq 0, ~x(t) \in {\bf R}^n , $$
\begin{equation}
x(\xi) = \varphi(\xi), ~if~~ \xi < 0. \label{e11}
\end{equation}

The equation (\ref{e11}) determines a differential operator ${\cal
M}$
\begin{equation}
 ({\cal M} x)(t) = \dot{x}(t) + \sum_{k=1}^r H_k (t) x[g_k(t)],
{}~x(\xi)=0, ~\xi < 0.
\label{e82}
\end{equation}

Suppose for this equation the hypotheses (a1)-(a4) hold.
By $C_{\cal M}$ we denote the Cauchy operator of this equation.

\begin{uess}
Suppose that for the operators ${\cal L}$ and ${\cal M}$
defined by (\ref{e4}) and (\ref{e82}) the following conditions
are satisfied.

1. The operators ${\cal L}$ and ${\cal M}$ act from $\depe$ to
$\lp$ and they are bounded.

2. $R({\cal M}) = \lp $, where $R({\cal M})$ is a range of values
of the operator ${\cal M} : \depe \rightarrow \lp $.

3. The operator ${\cal L} C_{\cal M} : \lp \rightarrow \lp$
is invertible.

Then $R({\cal L}) = \lp$ and $C$ acts from $\lp$ to $\depe$
and it is bounded.
\end{uess}

{\sl Proof.}
Consider an initial value problem
$$ {\cal L} x = f, ~x(0)=0, ~x(\tau_j)= B_j x(\tau_j -0), $$
where $f \in \lp$ is an arbitrary function.
Then $x = C_{\cal M} ({\cal L} C_{\cal M} )^{-1} f$
is the solution of this problem.
Therefore $x \in \depe$, hence $R({\cal L}) = \lp$.

Let $\depe ^0 = \{ x \in \depe : x(0)=0 \} $.
Then by the Banach theorem on an inverse operator
the operator $C: \lp \rightarrow \depe ^0$
is bounded.
So the operator $C : \lp \rightarrow \depe $
is also bounded.

Denote
\begin{equation}
\varphi^h (t) = \left\{ \begin{array}{ll}
\varphi [h(t)], & h(t) <0, \\
0, & h(t) \geq 0 , \end{array} \right.
g(t)= \sum_{k=1}^m A_k(t) \varphi^{h_k} (t).
\label{e12}
\end{equation}

\begin{guess}
Suppose the operators ${\cal L}$ and ${\cal M}$
defined by (\ref{e4}) and (\ref{e82})
satisfy the conditions of Lemma 5.

If the function $g$ defined by (\ref{e12}) is in $\lp$
then pair $(\lp,\depe)$ is admissible for the equation
(\ref{e1})-(\ref{e3}).

If there exists $\delta > 0$ such that
$t-h_k(t) < \delta$
and the restriction of $A_k$ to $[0, \delta ]$
belongs to $\lp [0, \delta], ~k=1, \dots, m,$
then the pair $(\lp,\depe)$ is admissible on the whole for
the equation (\ref{e1})-(\ref{e3}).
\end{guess}

{\sl Proof.}
Let $f \in \lp$ and $C$ be the Cauchy operator of
(\ref{e1})-(\ref{e3}).
By Lemma 1 solution $x$ of (\ref{e1})-(\ref{e3}) can be presented
as
\begin{equation}
x(t) = (Cf)(t) - (Cg)(t) + \sum_{0 \leq \tau_j \leq t} X(t,
\tau_j)
\alpha_j.   \label{e13}
\end{equation}

By Lemma 5 $Cf \in \depe, ~Cg \in \depe$.
Now we will establish $X(\cdot, \tau_j) \in \depe, ~j=1,2, \dots
.$
To this end denote
$$ Y_j (t) = X(t, \tau_j) - X_0 (t, \tau_j), $$
where $X_0 (t,s)$ is the fundamental matrix of (\ref{e6}),(\ref{e7})
and $a - I \ln b > 0$.
$$ \mbox{Let~~~} f_j (t) = - {\cal L} (X_0 (\cdot, \tau_j)) (t). $$

Then $Y_j$ is a solution of the problem
$${\cal L} y = f_j, ~ t \geq \tau_j,~ y(t) \in {\bf R}^{n \times
n}, $$
\begin{equation}
y(\tau_j)=0, ~ y(\tau_i) = B_i y(\tau_i - 0),
{}~~ i= j+1, \dots .
\label{e14}
\end{equation}

By Lemma 1 the solution of (\ref{e14}) can be presented as
$$ Y_j (t) = (Cf_j)(t), $$
hence
\begin{equation}
X(t,\tau_j) = X_0 (t, \tau_j) + (Cf_j)(t).
\label{e15}
\end{equation}

By Lemma 2  $X_0 (\cdot, \tau_j) \in \depe$.
Since by the hypothesis of the theorem pair $(\depe, \lp)$
is admissible for the operator ${\cal L}$ then $f_j \in \lp$.
Therefore by Lemma 5 $Cf_j \in \depe$.
Thus (\ref{e15}) implies $X(\cdot, \tau_j) \in \depe$
and (\ref{e13}) gives that a solution of (\ref{e1})-
(\ref{e3}) is in $\depe$.
Admissibility of the pair $(\lp, \depe)$ for the equation
(\ref{e1})-(\ref{e3}) is proven.

Suppose $t-h_k (t) < \delta$.
As $g$ is defined by (\ref{e12}) then $g(t)=0$ if $t> \delta$.
Since for $t \in [0, \delta]$ ~$A_k \in \lp [0, \delta] $
and $\varphi^{h_k} \in {\bf L}_{\infty} [0, \delta] $,
then $g \in \lp [0, \delta] $.
Therefore for $t \in [0, \infty)$ $g \in \lp [0, \infty) $.
Thus according to the above results the pair $(\lp, \depe)$
is admissible on the whole for (\ref{e1})-(\ref{e3}).
The proof of the theorem is complete.

\section{Admissibility and stability}

{}~~~~~This paper deals with exponential stability only.
Other types of stability and their connection with properties of
the fundamental matrix are presented in [14].

\underline{\sl Definition.}
The equation (\ref{e1})-(\ref{e3})
is said to be {\bf exponentially stable} if there exist positive
constants $N$ and $\lambda$ such that for any initial function
$\varphi, f=0$ and
$\alpha_1= \alpha_2 = \dots =0$
 for a solution $x$ of (\ref{e1})-(\ref{e3})
the inequality
$$ \parallel x(t) \parallel \leq
N e^{- \lambda t} \left( \sup_{t<0}
\parallel \varphi (t) \parallel + \parallel x(0) \parallel
\right)
$$
holds.

Thus the representation (\ref{e5}) yields the following
assertion (see [14]).

\begin{guess}
Suppose (a1)-(a6) hold and
there exist positive constants $N$ and $\lambda$
such that the fundamental matrix $X(t,s)$ satisfies the inequality
\begin{equation}
\parallel X(t,s) \parallel \leq N e^{- \lambda (t-s)},
{}~t \geq s > 0,
\label{e16}
\end{equation}
and there exists $\delta > 0$ such that
$t-h_k (t) < \delta, ~k=1, \dots ,m$.

Then equation (\ref{e1})-(\ref{e3}) is exponentially stable.
\end{guess}

The following theorem is a main result of this work.
It connects admissibility of the pair $(\lp, \depe)$
with stability of (\ref{e1})-(\ref{e3}).

\begin{guess}
Suppose for (\ref{e1})-(\ref{e3}) the hypotheses (a1)-(a6),
hold, $A_k \in \mepe , ~1 \leq p < \infty$, there exists
$\delta > 0$ such that $t-h_k (t) < \delta, ~k=1, \dots, m$
and for the initial function $\varphi \equiv 0$
the pair $(\lp, \depe)$ is admissible for this equation.

Then the equation (\ref{e1})-(\ref{e3}) is exponentially stable.
\end{guess}

{\sl Proof.}
By Theorem 3 it is sufficient to prove that the estimate
(\ref{e16}) exists.
In view of Lemma 1 the fundamental matrix $X(t,s)$ as a function of
$t$ for a fixed $s$ is a solution of the problem
$$\dot{x}(t) + \sum_{k=1}^m A_k (t) x[h_k(t)] = 0,
{}~t \geq s,~ x(t) \in {\bf R}^{n \times n},~ x(s)= E_n, $$
\begin{equation}
x(\xi)= 0, ~ \xi < s, ~ x(\tau_j) = B_j x(\tau_j - 0), ~ \tau_j >
s.
\label{e17}
\end{equation}

Denote
\begin{equation}
Y(t,s) = e^{\lambda(t-s)} X(t,s),
\label{e18}
\end{equation}
where $\lambda > 0$ is a certain number.
Thus
$$Y(s,s) = X(s,s) = E_n \mbox{~~and,~besides,~}
Y(\tau_j,s) = B_j Y(\tau_j - 0,s), ~\tau_j >s. $$

Denote
$$ {\cal L}_s x =
\dot{x} (t) + \sum_{k=1}^m A_k (t) x[h_k(t)], t \geq s,
{}~x(t) \in {\bf R}^{n \times n};$$ $$
{}~x(\xi)=0,~\xi<s,.$$
By substituting $x(t) = y(t) e^{- \lambda (t-s)}$
we obtain
$$({\cal L}_s x)(t) = e^{- \lambda (t-s)} \dot{y} (s) -
e^{- \lambda (t-s)} \lambda y(t) +
\sum_{k=1}^m e^{- \lambda [h_k (t) -s]} A_k (t) y[h_k(t)] =
$$ $$ =
e^{- \lambda (t-s)} \left\{ \dot{y} (t) +
\sum_{k=1}^m A_k (t) y[h_k(t)] + \right. $$
$$ \left. + \sum_{k=1}^m e^{\lambda [t- h_k (t)]}
A_k (t) y[h_k (t)] - \sum_{k=1}^m A_k (t) y[h_k(t)] - \lambda
y(t) \right\} = $$
$$ = e^{- \lambda (t-s)} \left\{
({\cal L}_s y)(t) - \lambda y(t) + \sum_{k=1}^m
\left[ e^{\lambda (t-h_k(t))} - 1 \right]
A_k (t) y[h_k(t)] \right\}. $$

Denote
$$({\cal T}_s y)(t) =  \sum_{k=1}^m \left[
e^{\lambda (t- h_k(t))} - 1 \right]
A_k (t) y[h_k(t)] - \lambda y (t), ~t \geq s, $$
$$ ({\cal M}_s y)(t) = ({\cal L}_s y)(t)
+ ({\cal T}_s y)(t). $$
Then
$$ ({\cal L}_s x)(t)  =
e^{- \lambda (t-s)} ({\cal M}_s y)(t) $$
and $Y(t,s)$ is a fundamental matrix of the problem
${\cal M}_0 y \! = \! 0, ~y(\tau_j) = B_j y(\tau_j - 0)$.

The corollary of Theorem 1 gives that the operator ${\cal L}_s$
acts from $\depe [s, \infty)$ to $\lp [s, \infty)$
and it is bounded.
By the hypothesis of the theorem a solution
of ${\cal L}_s x = f$ together with its derivative is in
$\lp [s, \infty )$ if provided $f \in \lp [s, \infty )$.
Therefore the Cauchy operator $C_s$ of this equation
acts from $\lp [s, \infty)$ to $\depe [s, \infty)$.

Denote $\depe ^0 [s, \infty) =
\{ x \in \depe [s, \infty) \mid x(s)=0 \} $.
By the hypotheses of the theorem the operator
${\cal L}_s : \depe ^0 [s, \infty) \rightarrow \lp [s, \infty)$
is bounded.
By Lemma 3 the space $\depe [s, \infty)$ is Banach,
therefore its closed subspace $\depe ^0 [s, \infty)$
is also Banach.
Thus by the Banach theorem on an inverse operator
the operator $C_s : \lp [s, \infty) \rightarrow \depe ^0
[s, \infty)$ and, consequently, $C_s : \lp [s, \infty)
\rightarrow \depe [s, \infty)$ is bounded.

By Theorem 1
$H_s^k$ acts from $\depe [s, \infty)$ to
$\lp [s, \infty)$, where $(H_s^k x)(t) = A_k (t) x(h_k (t)), ~x(\xi)=0, ~\xi
<s $.
{}From the assumption $t- h_k (t) <  \delta$
we obtain an estimate
$$ \parallel {\cal T}_s \parallel _{\depe [s, \infty )
\rightarrow  \lp [s, \infty)}
\leq \left( e^{\lambda \delta} - 1 \right)
\sum_{k=1}^m \parallel H_s^k \parallel_{\depe [s, \infty)
\rightarrow \lp [s, \infty )} + \lambda . $$

The operator ${\cal M}_s C_s = E + {\cal T}_s C_s$,
with $E$ being an identity operator,
has a bounded inverse operator in $\lp [s, \infty )$ if
\begin{equation}
\parallel {\cal T}_s C_s \parallel_{\lp [s, \infty)
\rightarrow \lp [s, \infty)} < 1.
\label{e19}
\end{equation}

We prove that for $\lambda$ being small enough (\ref{e19}) holds.
To this end
$$ \parallel \! {\cal T}_s C_s \! \parallel _{\lp \rightarrow \lp}
\leq \parallel \! {\cal T}_s \! \parallel_{\depe \rightarrow \lp}
\parallel \! C_s \! \parallel _{\lp \rightarrow \depe }
\leq \! \left[ (e^{\lambda \delta} - 1)
\sum_{k=1}^m \!\! \parallel \! H_s^k \! \parallel + \lambda \right]
\! \parallel \! C_s \! \parallel . $$
Therefore for $\lambda$ being small enough (\ref{e19}) holds,
where $\lambda$ is obviously independent of $s$ since
$\parallel H_s^k\parallel~
\leq~\parallel H_0^k \parallel,
\parallel~C_s~\parallel~\leq~\parallel~C~\parallel $.

Operators ${\cal L}_s$ and ${\cal T}_s$ act continuously
from $\depe [s, \infty)$ to $\lp [s, \infty )$.
Hence the operator ${\cal M}_s = {\cal L}_s  +  {\cal T}_s $
also possesses this property.
Thus by Lemma 5 the Cauchy operator $C_{\cal M}^s$
of the equation ${\cal M}_s y = f$
continuously acts from $\lp [s, \infty)$ to $\depe [s, \infty)$.

Similar to (\ref{e15}) we obtain
\begin{equation}
Y(t,s) = X_0 (t,s) + (C_{\cal M}^s  f_s )(t).
\label{e20}
\end{equation}
Here $f_s (t) = - {\cal M}_s (X_0 (\cdot , s)) (t), ~a - I \ln b
> 0 $.

Lemma 2 implies $X_0 (\cdot, s) \in \depe [s, \infty )$.
Moreover, this lemma gives the uniform estimate
$ \parallel f_s \parallel_{\lp [s, \infty)}
\leq K, $
with $K$ not depending on $s$.

Therefore we obtain estimates independent of $s$
$$ \parallel C_{\cal M}^s f_s  \parallel_{\depe [s, \infty)}
\leq K \parallel C_{\cal M} \parallel, $$
$$ \parallel C_{\cal M}^s f_s \parallel_{\lp [s, \infty)}
\leq K \parallel C_{\cal M} \parallel .$$
and
$$ \parallel \frac{d}{dt}
C_{\cal M}^s f_s \parallel _{\lp [s, \infty)}
\leq K \parallel C_{\cal M} \parallel. $$

Denote $z_s = C_{\cal M}^s f_s$.
Since $z_s (s)=0$, then
$z_s = C_0^s (\dot{z}_s + a z_s )$.
By Lemma~2 $C_0^s : \lp [s, \infty) \rightarrow
{\bf L}_{\infty} [s, \infty)$
is bounded, therefore
$$\parallel C_{\cal M}^s f_s \parallel_{{\bf L}_{\infty} [s,
\infty)} = \parallel z_s \parallel_{{\bf L}_{\infty} [s, \infty)}
\leq $$
$$ \leq \parallel C_0 \parallel_{\lp \rightarrow {\bf L}_{\infty}}
\left( \parallel \dot{z}_s \parallel_{\lp [s, \infty )}
+ a \parallel z_s \parallel _{\lp [s, \infty)} \right) \leq $$
$$ \leq (1+a) K \parallel C_0 \parallel ~\parallel C_{\cal M}
\parallel . $$
Hence the estimate of the norm of $C_{\cal M}^s f_s$
in ${\bf L}_{\infty} [s, \infty) $ does not depend on $s$.

By Lemma 2 and (\ref{e20}) there exists $N>0$
such that $$vraisup_{t,s > 0}
\parallel Y(t,s) \parallel \leq N < \infty . $$

Thus (\ref{e18}) implies the exponential estimate (\ref{e16})
for the fundamental matrix of (\ref{e1})-(\ref{e3}).
The proof of the theorem is complete.

\section{Explicit stability results}

We apply Theorems 2 and 4 to obtaining explicit conditions of
exponential stability and of existence of integrable solutions.
To this end we prove an auxiliary result.

\begin{uess}
Suppose there exist $\sigma >0$ and
$\rho > 0$ such that $\rho \leq \tau_{j+1} - \tau_j \leq \sigma,
{}~ \parallel B_j \parallel \leq B < 1 $.

Then for the fundamental matrix $X_1$ of the equation
\begin{equation}
\dot{x} (t) = f(t), ~ x(\tau_j) = B_j x(\tau_j - 0)
\label{e21}
\end{equation}
the inequality
\begin{equation}
\parallel X_1 (t,s) \parallel
\leq e^{- \eta(t -s - \sigma )}
\label{e22}
\end{equation}
holds, where $\eta  = - \frac{1}{\sigma} \ln B $.
\end{uess}

{\sl Proof.}
Under the hypotheses of the lemma (see [13])
$$\parallel X_1 (t,s) \parallel \leq \left\{
\begin{array}{ll}
e^{- \eta (t-s)}  & , t -s > \sigma, \\
1  & , 0 < t - s \leq \sigma . \end{array} \right. $$
This immediately yields (\ref{e22}).

\vspace{5 mm}

Denote
$$ A_k^{\eta} (t) = A_k (t) e^{\eta (t- h_k (t))} . $$

\begin{guess}
Suppose for the equation (\ref{e1})-(\ref{e3})
the hypotheses (a3),(a4) and

(b1) $ f \in {\bf L}_1 , ~A_k^{\eta} \in {\bf M}_1 $;

(b2) $0 < \rho \leq \tau_{j+1} - \tau_j \leq \sigma ;$

(b3) $ \parallel B_j \parallel \leq B < 1 $;

(b4) $g  \in {\bf L}_1, $
where $g$ is defined by (13);

(b5) $e^{\eta ( \sigma + 1) } \sum_{k=1}^m
\parallel A_k^{\eta} \parallel_{{\bf M}_1} \leq 1 - e^{- \eta}
$,
where $\eta = - \frac{1}{\sigma} \ln B$,

hold.

Then for any solution $x$ of (\ref{e1})-(\ref{e3})
$x \in {\bf L}_1, ~\dot{x} \in {\bf L}_1 $.

\end{guess}

\begin{guess}
Suppose for the equation (\ref{e1})-(\ref{e3})
the hypotheses (a3),(a4) and

(c1) $f \in {\bf L}_1, ~ A_k \in {\bf M}_1 $;

(c2) $0 < \rho \leq \tau_{j+1} - \tau_j \leq \sigma $;

(c3) $\parallel B_j \parallel \leq B < 1 $ ;

(c4) there exists $\delta > 0$ such that $t - h_k(t) < \delta $ ;

(c5) $ e^{\eta (\sigma + \delta + 1)}
\sum_{k=1}^m \parallel A_k \parallel_{{\bf M}_1}
\leq 1-e^{- \eta},$ where $ \eta =
- \frac{1}{\sigma} \ln B$ ,

hold.

Then the equation (\ref{e1}) - (\ref{e3}) is exponentially
stable.
\end{guess}

\underline{\sl Proof of Theorem 5.}
First we note that the hypotheses of the theorem imply (a1)-(a6).
In particular, (b2) implies (a1) and (a6).
By Theorem 1 the hypotheses of the theorem ensure admissibility
of the pair $({\bf D}_1, {\bf L}_1)$ for operator ${\cal L}$
defined by (4).

The hypotheses of Theorem 2 are satisfied if operator
${\cal L}C_{\cal M}:{\bf L}_1~\rightarrow{\bf L}_1$
is invertible, where
$C_{\cal M}$ is the Cauchy operator of the problem (\ref{e21}).

Evidently
${\cal L}C_{\cal M}=E+T$, where
$$(Tz)(t)= \sum_{k=1}^m A_k (t)
\int_0^{h_k^+(t)} X_1 (h_k(t),s)z(s)ds. $$

Lemma 6 gives that the operator $C_{\cal M}$
acts from ${\bf L}_1$ to ${\bf D}_1$.
Since by the hypothesis of the theorem $A_k^{\eta}\in~{\bf M}_1$,
then from the equality $T=H C_{\cal M},$
where $H$ is defined by (\ref{e9}), and from Theorem 1
the operator $T$ acts in ${\bf L}_1$.

Let estimate the norm of operator $T$:
$$ \parallel Tz \parallel_{{\bf L}_1}  \leq
\sum_{k=1}^m \int_0^{\infty}
\parallel A_k(t) \parallel \int_0^{h_k^+(t)}
e^{- \eta (h_k (t) -s- \sigma)}
\parallel z(s) \parallel ds~dt \leq $$
$$ \leq e^{\eta \sigma}
\sum_{k=1}^m \int_0^{\infty} \int_0^t
\parallel A_k (t) e^{\eta (t-h_k(t))}
e^{- \eta (t-s)}
\parallel z(s) \parallel ds~dt = $$
$$ = e^{\eta \sigma} \sum_{k=1}^m \int_0^{\infty}
\left( \int_s^{\infty} \parallel A_k^{\eta } (t) \parallel
e^{- \eta  t} dt \right) e^{\eta  s}
\parallel z(s) \parallel ds \leq $$
$$ \leq e^{\eta \sigma}
\sum_{k=1}^m \int_0^{\infty}
\left( \sum_{i=[s]}^{\infty}
\int_i^{i+1} \parallel A_k^{\eta} (t) \parallel e^{- \eta t}
dt \right) e^{\eta s} \parallel z(s) \parallel ds \leq $$
$$ \leq e^{\eta \sigma} \sum_{k=1}^m \parallel A_k^{\eta}
\parallel_{{\bf M}_1} \int_0^{\infty}
\sum_{i=[s]}^{\infty} e^{-\eta i} e^{\eta s}
\parallel z(s) \parallel ds = $$
$$ = e^{\eta \sigma} \sum_{k=1}^m
\parallel A_k^{\eta} \parallel_{{\bf M}_1}
\frac{e^{\eta}}{1-e^{- \eta}} \parallel z \parallel_{{\bf L}_1}. $$

The hypothesis (b5) implies
$\parallel T \parallel_{{\bf L}_1 \rightarrow {\bf L}_1}<1$,
therefore
${\cal L}C_{\cal M}:{\bf L}_1~\rightarrow~{\bf L}_1$
is invertible.
Hence all hypotheses of Theorem 2 hold.
The proof of the theorem is complete.
\vspace{5 mm}

\underline{\sl Proof of Theorem 6.}
The hypothesis (c4) implies $\varphi (h_k(t))=0$
for $t> \delta$.
Thus (b1),(b4),(b5) and other hypotheses of Theorem 5 hold.
By Theorem 4 the equation (\ref{e1})-(\ref{e3})
is exponentially stable.

\underline{\sl Example.}
Consider a scalar equation
$$\dot{x}(t)+a(t)x(\lambda t)=f(t),~t \geq 0,~ 0 < \lambda < 1,$$
\begin{equation}
 x(j)=bx(j-0),~j=1,2, \dots,~ \mid b \mid <1 .
\end{equation}

Since $h(t)= \lambda t \geq 0$ then one may assume $\varphi
\equiv 0$.

The constant $\eta$ defined in (b5)
is $\eta  = - \ln b$.
Therefore by Theorem 5 all solutions of (24)
are in ${\bf L}_1 $ for any $f \in {\bf L}_1,$
 i.e. they are integrable on the half-line
if
$$
a^{\eta} (t) = a(t) e^ {[(\lambda - 1) \ln b] t}  \in {\bf M}_1
\mbox{~~ and ~~}
\parallel a^{\eta}  \parallel _{{\bf M}_1} \leq (1-b)b^2 . $$

\section{Proofs of Lemmas 3 and 4}

{\bf Lemma 3.}
{\it Suppose (a5) and (a6) hold.
Then $\depe, ~1 \leq p < \infty $,
is a Banach space. }

{\sl Proof.}
Let $\{x_j \}$ be a fundamental sequence in $\depe$, i.e.
$$\lim_{k,i \rightarrow \infty}
\parallel x_k - x_i \parallel _{\depe} = 0. $$

First we will prove that $\{ x_j (0)\}$
converges in ${\bf R}^n$.

The convergence $\parallel y_j \parallel_{\depe [0, \infty)}
\rightarrow 0$ implies
$\parallel y_j \parallel_{\depe [0,t_0]} \rightarrow 0$ for any
$t_0 >0$.
Hence
$\parallel y_j \parallel_{{\bf L}_1 [0,t_0]} \rightarrow 0,
{}~ \parallel \dot{y}_j \parallel_{{\bf L}_1 [0,t_0]} \rightarrow 0$.
Therefore for $t_0 < \tau_1$  and $y_j = x_k -x_i$ we have
$$ \lim_{k,i \rightarrow \infty}
\parallel x_k - x_i \parallel_{{\bf L}_1 [0,t_0]} = 0,
\lim_{k,i \rightarrow \infty}
\parallel \dot{x}_k - \dot{x}_i \parallel_{{\bf L}_1 [0,t_0]} = 0 .
$$
Consider an identity
$$ x_k (t) - x_i (t) = x_k (0) - x_i (0) +
\int_0^t [\dot{x}_k (s) - \dot{x}_i (s)] ds. $$
Since
$$\lim_{k,i \rightarrow \infty}
\parallel x_k - x_i \parallel_{{\bf L}_1 [0,t_0]} = 0 \mbox{ and}$$
$$\lim_{k,i \rightarrow 0} \int_0^{t_0} \int_0^t
\mid \dot{x}_k (s) - \dot{x}_i (s) \mid ~ ds~dt \leq
t_0 \lim_{k,i \rightarrow \infty}
\parallel \dot{x}_k - \dot{x}_i \parallel_{{\bf L}_1 [0,t_0]} = 0,
$$
then
$$\lim_{k,i \rightarrow \infty} \parallel x_k (0) - x_i (0)
\parallel_{{\bf L}_1 [0,t_0]} = 0. $$
Hence
$$\lim_{k,i \rightarrow \infty} \parallel x_k (0) - x_i (0)
\parallel_{{\bf L}_1 [0,t_0]} =
t_0 \lim_{k,i \rightarrow \infty}
\parallel x_k (0) - x_i (0) \parallel_{{\bf R}^n} = 0, $$
i.e. the sequence $\{x_j (0)\}$ is fundamental in ${\bf R}^n$.
Therefore there exists $\beta \in {\bf R}^n$ such that
$\lim_{j \rightarrow \infty} x_j (0) = \beta $ .

Let $f_j = {\cal L}_0 x_j$, where operator ${\cal L}_0$
is defined by (\ref{e6}), $\nu = a - I \ln b > 0$.
Then by Lemma 1
\begin{equation}
x_j (t) = X_0 (t,0) x_j (0) + \int_0^t X_0 (t,s) f_j (s) ds.
\label{e80}
\end{equation}
Lemma 2 yields
$$ \parallel X_0 (t,0) \parallel \leq e^{- \nu t} . $$
Since $\dot{X}_0 (t,0) + a X_0 (t,0)= 0$ then
$$ \parallel \dot{X}_0 (t,0) \parallel \leq a e^{- \nu t} .$$
Therefore the sequence $\{X_0 (t,0) x_j (0)\}$ converges in
$\depe$ to the function $X_0 (t,0) \beta $.

By Lemma 2 we obtain that the operators ${\cal L}_0 : \depe
\rightarrow \lp$ and $C_0 : \lp \rightarrow \depe $ are bounded.
To this end denoting $x = C_0 f$ we obtain
$$ \parallel {\cal L}_0 x \parallel_{\lp} \leq
\parallel \dot{x} \parallel_{\lp} + a \parallel x \parallel_{\lp}
\leq (1+a) \parallel x \parallel _{\depe}. $$
By Lemma 2 operator $C_0$ is bounded in $\lp$ [1], hence
$$ \parallel C_0 f \parallel_{\depe} = \parallel x
\parallel_{\lp} + \parallel \dot{x} \parallel_{\lp} \leq
\parallel C_0 \parallel
_{\lp \rightarrow \lp}
 \parallel f \parallel _{\lp} +
\parallel f-ax \parallel _{\lp} \leq $$
$$ \leq
[ 1+ \parallel C_0 \parallel (1+a)] \parallel f \parallel_{\lp}. $$

Since ${\cal L}_0 : \depe \rightarrow \lp$
is continuous and ${\cal L}_0 x_j = f_j$ then
$\{ f_j \}$ is a fundamental sequence.
Therefore there exists $f \in \lp$ such that $\lim_{j \rightarrow
\infty} f_j = f$.

Let $\tilde{x} = C_0 f,~ \tilde{x}_j = C_0 f_j $.
The continuity of the operator $C_0: \lp \rightarrow \depe$
implies
$ \parallel \tilde{x}_j - \tilde{x} \parallel_{\depe} \rightarrow 0. $

{}From here sequence
$$ x_j (t) = X_0 (t,0) x_j (0) + \tilde{x}_j (t) $$
converges in $\depe$ to
$$ x (t) = X_0 (t,0) \beta + \tilde{x} (t). $$
The proof of the lemma is complete.
\vspace{5 mm}

\underline{\sl Lemma 4.}
{\it Let } $a> I \ln b$.

{\it Then the set}
$$\tilde{\depe} = \{ x \in {\bf PAC} \mid  \dot{x} + ax \in \lp ,
{}~x(\tau_j) = B_j x(\tau_j -0) \} $$
{\it coincides with $\depe$.
Besides the norm }
\begin{equation}
\parallel x \parallel_{\tilde{\depe}} = \parallel x(0) \parallel +
\parallel \dot{x} + ax \parallel_{\lp}
\label{e25}
\end{equation}
{\it is equivalent to the norm in } $\depe$.

{\sl Proof.}
Let $x \in \tilde{\depe}$ and $z = \dot{x} + ax$.
Then $x(t) = X_0 (t,0) x(0) + (C_0 z )(t)$.
By Lemma 2  $z \in \lp$
implies $x \in \lp$.
Hence $\dot{x}~=z-ax~\in~\lp,$ thus $x \in \depe$.

Let $x \in \depe$.
Then the inequality
$$ \parallel \dot{x} + ax \parallel_{\lp}
\leq (1+a) \parallel x \parallel_{\depe} $$
implies
$\dot{x} + ax \in \lp$.
Hence $x \in \tilde{\depe}$.
Thus $\tilde{\depe}=\depe$.

Formula (\ref{e25}) defines a norm in $\depe$.
In fact if $\parallel x\parallel_{\tilde{\depe}}~=0$
then $\dot{x}+ax=0,~ x(0)=0$.
Then by Lemma 1 on uniqueness of a solution $x=0$.

Let us prove that the space $\depe$
endowed with the norm $\parallel \cdot \parallel _{\tilde{\depe}}$
is complete.
Suppose $\{ x_j \}$ is a fundamental sequence
by this norm.
Denote $y_j = \dot{x}_j+ax_j$.
Then the convergence
$$\parallel x_k (0) - x_i (0) \parallel +
\parallel y_k - y_i \parallel _{\lp} \rightarrow 0
\mbox{~~for~~} k,i \rightarrow \infty$$
implies $\{x_j (0)\}$ is fundamental in ${\bf R}^n$
and $\{y_j \}$ is fundamental in $\lp$.
Therefore these sequences converge in the corresponding spaces.

Consider the equality
\begin{equation}
x_j (t) = X_0 (t,0) x_j (0) + (C_0 y_j)(t).
\label{e60}
\end{equation}
We will prove that the operator $C_0: \lp \rightarrow
\tilde{\depe} $ is bounded.
Let $x=C_0 f$.
Then $x(0)=0$ and
$$ \parallel C_0 f \parallel_{\tilde{\depe}} =
\parallel \dot{x} + ax \parallel_{\lp} =
\parallel f \parallel_{\lp} . $$
Boundedness of $C_0 : \lp \rightarrow \tilde{\depe}$
and the equality (\ref{e60}) yield the convergence of $\{x_j \}$
in $\tilde{\depe}$.
Consequently this space is complete.

Consider sets
$$\depe^0 = \{ x \in \depe \mid x(0)=0 \} , $$
$$U_n = \{ x=X_0 (t,0) \alpha \mid ~ \alpha \in {\bf R}^n \} . $$
The space $U_n$ is $n$-dimensional, isomorphic to ${\bf R}^n$
and $U_n \subset \depe$.
Since
$$x(t) = X(t,0)x(0) + \int_0^t X_0 (t,s) [\dot{x} (s) + ax(s)] ds
, $$
then $\depe$ is aljebraically isomorphic to the direct sum
$\depe^0 \oplus U_n$.

Since $U_n$ is finite-dimensional then [15]
the subspace $\depe^0$ is closed in $\depe$ and in $\tilde{\depe}$.

First we will prove equivalence of norms
$\parallel \cdot \parallel_{\depe}$ and $\parallel
\cdot \parallel_{\tilde{\depe}}$ in $\depe^0$.
Let $x~\in~\depe^0$.
To this end
$$\parallel x \parallel_{\tilde{\depe}} =
\parallel \dot{x} + ax \parallel_{\lp} \leq (1+a)
\parallel x \parallel_{\depe}. $$

{}From here and from the fact $D_p^0$ is a Banach space with
both norms we obtain [15] that in $D_p^0$ these norms
are equivalent.

Let $P_1$ and $P_2$ be projectors to subspaces $\depe^0$
and $U_n$ correspondingly.
$\depe^0$ is closed, therefore these projectors
are bounded operators in $\depe$ and $\tilde{\depe}$.
Let $\parallel x_j \parallel_{\depe} \rightarrow 0$.
Then the relations
$$x_j=P_1x_j+P_2x_j,
 ~ \parallel P_i x_j \parallel _{\depe} \leq
 \parallel P_i \parallel \parallel x_j \parallel _{\depe},
 ~i=1,2, $$
 imply $\parallel P_i x_j \parallel_{\depe} \rightarrow 0,
 ~i=1,2$.
As $P_1 x_j \in \depe^0,$
and in $\depe^0$ the norms
$\parallel \cdot \parallel_{\depe}$ and
$\parallel \cdot \parallel_{\tilde{\depe}}$
are equivalent, then
$\parallel P_1 x_j \parallel_{\tilde{\depe}} \rightarrow 0$.

Besides this $P_2 x_j \in U_n$.
The space $U_n$ is finite-dimensional and all the norms in it
are equivalent.
Thus
$\parallel P_2 x_j \parallel_{\tilde{\depe}} \rightarrow 0$.
Consequently,
$$ \parallel x_j \parallel_{\tilde{\depe}} \leq
\parallel P_1 x_j \parallel_{\tilde{\depe}}
+ \parallel P_2 x_j \parallel_{\tilde{\depe}} \rightarrow 0 .$$
Therefore the norms $\parallel \cdot \parallel_{\depe}$
and $\parallel \cdot \parallel_{\tilde{\depe}}$ are equivalent,
which completes the proof.

\end{document}